\documentclass[12pt]{article}
\usepackage{epsfig}
\usepackage[usenames]{color}
\usepackage{verbatim}
\usepackage{hyperref}
\usepackage{multicol}
\usepackage{comment}
\usepackage{float}
\usepackage{graphicx}
\usepackage[utf8]{inputenc}

\usepackage{amsfonts,amsmath}

\def \beq {\begin{eqnarray}}
\def \eeq {\end{eqnarray}}
\def \beqn {\begin{eqnarray*}}
\def \eeqn {\end{eqnarray*}}

\newcommand{\halmos}{\rule{1ex}{1.4ex}}
\newcounter{for}[section]

\newtheorem{itlemma}{Lemma}[section]
\newtheorem{itproposition}[itlemma]{Proposition}
\newtheorem{theorem}[itlemma]{Theorem}
\newtheorem{itcorollary}[itlemma]{Corollary}
\newtheorem{itremark}[itlemma]{Remark}
\newtheorem{itremarks}[itlemma]{Remarks}
\newtheorem{itdefinition}[itlemma]{Definition}
\newtheorem{itexample}[itlemma]{Example}

\newenvironment{fact}{\begin{itfact}\rm}{\end{itfact}}
\newenvironment{claim}{\begin{itclaim}\rm}{\end{itclaim}}
\newenvironment{lemma}{\begin{itlemma}}{\end{itlemma}}
\newenvironment{remark}{\begin{itremark}\rm}{\end{itremark}}
\newenvironment{remarks}{\begin{itremarks} \rm}{\end{itremarks}}
\newenvironment{corollary}{\begin{itcorollary}}{\end{itcorollary}}
\newenvironment{proposition}{\begin{itproposition}}{\end{itproposition}}
\newenvironment{definition}{\begin{itdefinition}\rm}{\end{itdefinition}}
\newenvironment{example}{\begin{itexample}\rm}{\end{itexample}}
\newenvironment{proof}{\noindent {\em Proof}.\ \
}{\hspace*{\fill}$\halmos$\medskip}
\newcommand{\be}[1]{\addtocounter{for}{1} \begin{equation}\label{#1}}
\newcommand{\ee}{\end{equation}}
\newcommand{\bl}[1]{\begin{lemma}\label{#1}}
\newcommand{\br}[1]{\begin{remark}\label{#1}}
\newcommand{\brs}[1]{\begin{remarks}\label{#1}}
\newcommand{\bt}[1]{\begin{theorem}\label{#1}}
\newcommand{\bd}[1]{\begin{definition}\label{#1}}
\newcommand{\bp}[1]{\begin{proposition}\label{#1}}
\newcommand{\bc}[1]{\begin{corollary}\label{#1}}
\newcommand{\bfact}[1]{\begin{fact}\label{#1}}
\newcommand{\bex}[1]{\begin{example}\label{#1}}
\newcommand{\ec}{\end{corollary}}
\newcommand{\efact}{\end{fact}}
\newcommand{\eex}{\end{example}}
\newcommand{\el}{\end{lemma}}
\newcommand{\er}{\end{remark}}
\newcommand{\ers}{\end{remarks}}
\newcommand{\et}{\end{theorem}}
\newcommand{\ed}{\end{definition}}
\newcommand{\ep}{\end{proposition}}
\newcommand{\epr}{\end{proof}}
\newcommand{\bpr}{\begin{proof}}
\newcommand{\bcl}[1]{\begin{claim}\label{#1}}
\newcommand{\ecl}{\end{claim}}

\newcommand{\ecs}{\end{corollary}}
\newcommand{\eers}{\end{exercise}}
\newcommand{\eexs}{\end{example}}
\newcommand{\eems}{\end{example}}
\newcommand{\els}{\end{lemma}}
\newcommand{\eles}{\end{lemmaex}}
\newcommand{\ets}{\end{theorem}}
\newcommand{\eds}{\end{definition}}
\newcommand{\eps}{\end{proposition}}
\newcommand{\bi}{\begin{itemize}}
\newcommand{\ei}{\end{itemize}}
\newcommand{\ben}{\begin{enumerate}}
\newcommand{\een}{\end{enumerate}}

\def\vbar{\mathchoice{\vrule height6.3ptdepth-.5ptwidth.8pt\kern-.8pt}
   {\vrule height6.3ptdepth-.5ptwidth.8pt\kern-.8pt}
   {\vrule height4.1ptdepth-.35ptwidth.6pt\kern-.6pt}
   {\vrule height3.1ptdepth-.25ptwidth.5pt\kern-.5pt}}
\def\fudge{\mathchoice{}{}{\mkern.5mu}{\mkern.8mu}}
\def\bbc#1#2{{\rm \mkern#2mu\vbar\mkern-#2mu#1}}
\def\bbb#1{{\rm I\mkern-3.5mu #1}}
\def\bba#1#2{{\rm #1\mkern-#2mu\fudge #1}}
\def\bb#1{{\count4=`#1 \advance\count4by-64 \ifcase\count4\or\bba A{11.5}\or
   \bbb B\or\bbc C{5}\or\bbb D\or\bbb E\or\bbb F \or\bbc G{5}\or\bbb H\or
   \bbb I\or\bbc J{3}\or\bbb K\or\bbb L \or\bbb M\or\bbb N\or\bbc O{5} \or
   \bbb P\or\bbc Q{5}\or\bbb R\or\bbc S{4.2}\or\bba T{10.5}\or\bbc U{5}\or
   \bba V{12}\or\bba W{16.5}\or\bba X{11}\or\bba Y{11.7}\or\bba Z{7.5}\fi}}

\def \S {{\cal{S}}}

\newcommand{\ba}[1]{\addtocounter{for}{1} \begin{eqnarray}\label{#1}}
\newcommand{\ea}{\end{eqnarray}}

\def\sqr#1#2{{\vcenter{\vbox{\hrule height .#2pt
                             \hbox{\vrule width .#2pt height#1pt \kern#1pt
                                   \vrule width .#2pt}
                             \hrule height .#2pt}}}}

\def\pmb#1{\setbox0=\hbox{#1}%
   \kern-.025em\copy0\kern-\wd0
   \kern.05em\copy0\kern-\wd0
   \kern-.025em\raise.0433em\box0 }
\def\sqr#1#2{{\vcenter{\vbox{\hrule height.#2pt
     \hbox{\vrule width.#2pt height#1pt \kern#1pt
   \vrule width.#2pt}\hrule height.#2pt}}}}

\def\N{{\mathbb N}}
\def\Z{{\mathbb Z}}
\def\R{{\mathbb R}}

\def\bs{\backslash}

\def\bs{\backslash}
\def\reff#1{(\ref{#1})}

\newcommand {\acc}[1] {\left\{ {#1} \right\}}
\newcommand {\pare}[1] {\left( {#1} \right)}

\newcommand{\1}{{\text{\Large $\mathfrak 1$}}}

\def\cV{\mathcal{V}}
\def\cS{\mathcal{S}}
\def\cN{\mathcal{N}}

\begin{document}
\title{\bf On outer fluctuations for internal DLA}

\author{Amine Asselah \thanks{
Universit\'e Paris-Est, LAMA (UMR 8050), UPEC, UPEMLV, CNRS, F-94010, Cr\'eteil, France; amine.asselah@u-pec.fr} \and
Alexandre Gaudilli\`ere \thanks{Aix-Marseille Universit\'e, CNRS, Centrale Marseille, I2M, UMR 7373, 13453 Marseille, France;  
alexandre.gaudilliere@math.cnrs.fr}}
\date{}
\maketitle
\begin{abstract}
We had established in \cite{AG} inner and outer fluctuation
for the internal DLA cluster when all walks are launched from
the origin. In obtaining the outer fluctuation,
we had used a deep lemma of Jerison, Levine and Sheffield of \cite{JLS},
which estimate roughly the possibility of {\it fingering},
and had provided in \cite{AG} a simple proof using an interesting estimate
for crossing probability for a simple random walk. The application
of the crossing probability to the fingering for the internal DLA cluster
contains a flaw discovered recently, that we correct in this note. 
We take the opportunity to make a self-contained exposition.
\end{abstract}
\section{Introduction}
In this short note, we correct a mistake in
an alternative proof we gave in \cite{AG} (proof of Lemma 1.5) of
Lemma A of David Jerison, Lionel Levine
and Scott Sheffield in \cite{JLS}. This
result bounds the probability the cluster of internal DLA, with particles 
starting outside a ball, eventually covers the center of this ball
when the number of particles is small compared to the volume of the ball.
Lemma A controls the possibility the cluster makes fingers protruding
out of the spherical shape it likely adopts.
This estimate is used in turn to produce an
outer error bound when all internal DLA particles are 
launched from the origin of $\Z^d$ and $d\ge 2$.

Recently Lionel Levine and Yuval Peres noticed a flaw in our (simple)
proof. We correct this flaw by adding one step to our initial proof, 
and since we believe the result is of independent interest, 
we take the opportunity to present 
a self-contained argument by including an estimate of the
probability one random walk crosses a shell while
staying inside a given region in term of its volume. 
This estimate is Lemma 1.6 of \cite{AG} and concerns one
random walk.

The walk is denoted $S: \N\to\Z^d$, and we call 
$P_z$ the law of the simple random walk when $S(0)=z$.
If $\Lambda$ is a subset of $\Z^d$, $T(\Lambda)$ denotes
the hitting time of $\Lambda$. Also, we call $\|\cdot\|$ the
euclidean distance, $B(z,\rho)$
the trace on $\Z^d$ of the
ball of radius $\rho>0$ and center $z\in \R^d$, and 
$\partial B(z,\rho)$ the boundary of $B(z,\rho)$, that is
$\partial B(z,\rho):=\{y\in \Z^d\bs B(z,\rho):\ \exists x\in B(z,\rho),
\ \|x-y\|=1\}$. We can now state one key Lemma of \cite{AG}.
\bl{lem-1}
Assume dimension is $d\ge 2$. There are constants
$\kappa_d>0,C_d>1$ such that for any $r>0$, $0<h<r/2$
and $z\in \partial B(0,r)$,
and $\cV\subset B(0,r)\bs B(0,r-h)$ we have
\be{ineq-main}
P_z\big(T(B(0,r-h)) < T(\cV^c)\big)\le 
C_d\exp\pare{-\kappa_d
\pare{\frac{h^d}{|\cV|}}^{\frac{1}{d-1}}}
\ee
\el
Note that the estimate \reff{ineq-main} is useful
only if $h^d/|\cV|$ is large enough to counter $C_d$.
To state our next result, we need to introduce internal DLA and
more notation.
A configuration of random walkers is denoted $\eta$ and is an
element of $\N^{\Z^d}$. The number of walks in $\eta$ is $|\eta|=
\sum_{z\in \Z^d} \eta(z)$. The cluster of internal DLA made of
random walks initially in $\eta$ with $|\eta|<\infty$, is 
itself a configuration of $\{0,1\}^{\Z^d}$ 
that we build inductively
as follows. Choose an arbitrary ordering of the random walks,
and run the walks one at a time following their order. 
When the running walk steps on
an empty site, it stops, or {\it settles} and the site 
becomes part of the cluster. At this moment, 
send the next random walk until it {\it settles} and so on. 
This cluster has a law independent of the ordering of the 
walks. This is the celebrated {\it abelian property}. 
Note that only one random walk settles in
each site of the final cluster that we call $A(\eta)$ and
which can be seen as a subset of $\Z^d$. The random walks with the rules
for settling are called {\it explorers}.

We can now state
a weaker result than Lemma 1.5 of \cite{AG}.
It is a direct consequence of the (corrected) proof of Lemma 1.5.

\bl{lem-levine} Assume that $d=2$.  
There are $a,\kappa_2>0$, such that for $r$ large enough, 
and for any configuration $\eta$ of explorers outside of $B(0,2r)$,
with $|\eta|\le a r^2/\log^2 r$, we have
\be{ineq-low2}
P\big(0\in A(\eta)\big)\le \exp\big(-\kappa_2 \frac{r^2}{\log(r)}\big).
\ee
Assume $d\ge 3$. 
There are $a,\kappa_d>0$, such that for $r$  large enough, and
for any configuration $\eta$ of explorers outside of $B(0,2r)$,
with $|\eta|\le a r^2/\log r$, we have
\be{ineq-low3}
P\big(0\in A(\eta)\big)\le \exp\big(-\kappa_d r^2\big).
\ee
\el
Finally, we state Lemma A of Jerison, Levine and Sheffield,
(and Lemma 1.5 of \cite{AG}) as a Corollary of Lemma~\ref{lem-levine}.

\bc{cor-19} The inequalities \reff{ineq-low2} and \reff{ineq-low3}
hold as soon as $r$ is large enough and any configuration
$\eta$ of explorers outside $B(0,2r)$ with
$|\eta|\le \epsilon r^d$, for some positive constant $\epsilon$
depending only on dimension.
\ec

We prove Lemma~\ref{lem-1} in Section~\ref{sec-one}.
We prove Lemma~\ref{lem-levine} in Section~\ref{sec-many},
and we correct our proof in Section~\ref{sec-old}.

\section{One explorer crossing a shell}\label{sec-one}
In this Section, 
we reproduce the short proof of Lemma~\ref{lem-1} of \cite{AG}.

Take a positive integer $n<h$, and consider a partition
of the shell $\S:=B(0,r)\bs B(0,r-h)$ into $n$ shells $\{\S_k,\ k<r\}$
of width $2\delta:=h/n$. For $k<n$, set
$\Sigma_k:=\partial B(0,r-(2k+1)\delta)$. Let $\{S(n),n\in \N\}$
be the underlying random walk with which we build an explorer.

With each explorer of internal DLA is associated a so-called
{\it flashing explorer} which can settle only on some random sites,
when they are empty.
Thus, we define the random sites $\{Z_k, 0\leq k<n\}$ as follows.
For each $k<n$, we draw a continuous random variable $R_k$ on $[0,\delta]$
with density in $x\in [0,\delta]\mapsto {dx^{d-1}}/{\delta^d}$, and
$Z_k$ is the exit site of $S$ from $B(S(T(\Sigma_k)),R_k)$
after time $T(\Sigma_k)$. Then, the flashing
explorer settles on the first $Z_k$ not belonging to $\cV$.
The purpose of the flashing construction is that
(i) the flashing sites are each distributed almost uniformly
inside the ball $B(S(T(\Sigma_k)),\delta)$ 
(and this is Proposition 3.1 of \cite{AG0}), and
(ii) $P_z(T(B(0,r)) < T(\cV^c))$ is bounded above
by the probability that the flashing explorer crosses $\S$.

Now, for a small $\beta$ to be chosen later, we say that $y\in \Sigma_k$
has a {\it dense neighborhood} if $|B(y,\delta) \cap \cV| 
\geq \beta \delta^d$, and we call $D_k\subset \Sigma_k$ the set
of such $y$. 
There is $\kappa>0$ such that knowing that the explorer
has crossed $D_1,\dots,D_{k-1}$, we have the following.
\begin{itemize}
\item If $S(T(\Sigma_k))\not\in D_k$, then the probability that
the explorer does not settle in $\S_k$ is smaller than $\kappa\beta$
(Proposition 3.1 of \cite{AG0}).
\item The probability that $S(T(\Sigma_k))\in D_k$ is smaller than
$\kappa |D_k|/h^{d-1}$ (see Lemma 5 of \cite{LBG})
uniformly over $Z_{k-1}$.
\end{itemize}
If the explorer has crossed $\S$, the flashing has also done so, 
which means that $Z_k\in \cV$ for all $k<n$.
By successive conditioning, we obtain
\be{alex-1}
P_z\pare{T(B(0,r-h)) < T(\cV^c)} \leq \prod_{k<n}
\pare{\kappa\beta+\frac{\kappa |D_k|}{\delta^{d-1}}}.
\ee
By the arithmetic-geometric inequality, we obtain
\be{arith-geo}
P_z\pare{T(B(0,r-h)) < T(\cV^c)} \leq
\left(\kappa\beta + \frac{\kappa}{n}
\sum_{k<n}\frac{|D_k|}{\delta^{d-1}}\right)^{n}.
\ee
Note that while each $y\in D_k$ satisfies $|B(y,\delta)\cap \cV|
\ge \beta \delta^d$, each site in $B(y,\delta)\cap \cV$ 
is at a distance less than $\delta$
from a number of sites of $D_k$ of 
order at most $\delta^{d-1}$. Thus, for some $c>0$
\[
\sum_{k<n}\frac{\beta |D_k| \delta^d}{\delta^{d-1}} \leq c|\cV|,
\quad\mbox{i.e.,}\quad
\frac{1}{n}\sum_{k<n}\frac{|D_k|}{\delta^{d-1}}
\leq \frac{2c|\cV|}{\beta h \delta^{d-1}}
\qquad\big(\text{we recall}\quad \frac{1}{n\delta}=\frac{2}{h}\big).
\]
We choose now $\beta$ such that $4\kappa \beta<1$, 
and we choose the smallest $\delta$ such that 
\be{alex-3}
\delta^{d-1}\geq \frac{2c|\cV|}{\beta^2 h}.
\ee
Thus, \reff{arith-geo} reads
\be{alex-7}
P_z(T(B(0,r-h)) < T(\cV^c)) \leq \pare{\frac{1}{2}}^{h/(2\delta)}.
\ee
Requiring that $2\delta<h$ adds a constraint on $|\cV|$:
\be{alex-4}
|\cV|\le \frac{\beta^2}{2^dc} h^d.
\ee
Instead of including \reff{alex-4} as a condition of our Lemma, 
we find it more convenient to note that the probability 
we estimate is less than 1, so that we obtain \reff{ineq-main},
with constant $C_d$.

\section{Cloud of Explorers Crossing a Shell}\label{sec-many}
We prove in this Section Lemma~\ref{lem-levine}.
Our initial problem consists in estimating the
crossing probability of one explorer out of $\eta$
when the positions of $\eta$ lay in the boundary of a 
ball of radius $2r$. The key idea of the proof is 
to divide the shell $B(0,2r)\bs B(0,r)$
into a sequence of smaller shells with random widths.
The choice is however different in the proofs of Corollary~\ref{cor-19}
and of Lemma~\ref{lem-levine}.
In proving Corollary~\ref{cor-19}, the width of a shell
depends on the number of explorers that reach the shell boundary.
This is natural in view of Lemma~\ref{lem-1}:
to estimate the crossing probability of a shell, 
the quantity one needs to control
is the number of explorers having 
settled in this very shell, which is bounded by the number
of explorers arriving on its external boundary.
However, in proving Lemma~\ref{lem-levine}
the width depends on the number of explorers
having settled in the {\em previous} shell.

Also, instead of considering one configuration $\eta$,
 we find it
convenient to take a Poisson cloud of explorers. This 
simplifies some large deviation  estimates.

\subsection{Poisson Cloud}
Let $\prec$ be the usual partial order on the space of configurations,
and consider an
increasing sequence of configurations $(\zeta_k,\ k\in \N)$
on the boundary of $B(2r)$. We require also that this sequence satisfies
$|\zeta_k|=k$, and
\[
\forall k\ge |\eta|,\quad \eta\prec \zeta_k.
\]
We call $u(\eta)$ the probability that one explorer settles in 
$B(r)$ when starting with an initial configuration $\eta$.
Note that the event we consider is increasing: with more explorers,
it is easier to make one of them cross. In other words,
\[
\forall k\ge |\eta|,\quad u(\zeta_k)\ge u(\eta).
\]
Let $K$ be a Poisson variable of parameter $\lambda_0:=\epsilon r^d$. The sequence $(\zeta_k,\ k\in \N)$ being given, we have
\[
E[u(\zeta_K)]=\sum_{k\in \N} u(\zeta_k)\cdot \frac{e^{-\lambda_0}
\lambda_0^k}{k!}.
\]
Since $P(K>\frac{1}{2}\lambda_0)\ge 1/8$ when $\lambda_0\ge 1$,
and $k\mapsto u(\zeta_k)$ is increasing, we have that if $k^*$ is
the integer part of $\frac{1}{2}\lambda_0$, and $r$ is large
enough so that $k^*>1$, then
\[
u\big(\zeta_{k^*}\big)\cdot P(K>k^*)\le E[u(\zeta_K)]
\Longrightarrow u\big(\zeta_{k^*}\big)\le 8 E[u(\zeta_K)].
\]
Thus, we need to estimate the expectation $E[u(\zeta_K)]$, where
the number of explorers is Poisson, whereas the initial configuration
is arbitrary. 
To this end we will make a repeated use
of the following deviation bound for 
a generic Poisson random variable $K$
of parameter $\lambda$: for any positive $\theta$
\be{ld-poisson}
P(K > \theta) \geq \exp\left\{ -\theta \log\left(
\theta \over e\lambda \right) \right\}
\ee
Let $\eta_0:=\zeta_K$,
and subdivide the shell $B(2r)\bs B(r)$ into successive shells 
of random width $H_0,H_1,\dots$
defined by induction as follows. For some constant $\gamma$,
\[
H_0^d=\gamma |\eta_0|=\gamma K,\quad\text{and}\quad
\cS_0:=B(2r)\bs B(2r-H_0).
\]
Imagine we have labelled the explorers, and have send $k-1$ of them,
that we stop either if they settle or when they enter $B(2r-H_0)$.
Let $\cV\subset \S_0$ be the domain where some have settled,
whereas at most $k-1-|\cV|$ are stopped on entering $B(2r-H_0)$. 
The probability that the $k$-th explorer,
with $k\le K$, crosses $\cS_0$, knowing $K$ and $\cV$, 
is bounded as we use \reff{ineq-main} and $|\cV|\le K$,
\be{explorer-1}
\begin{split}
P\big(\text{ the $k$-th explorer crosses }\cS_0\
 |K,\cV\big)\le&
\sup_{z\in \partial B(2r)}\sup_{\cV \subset\cS_0:\ |\cV|\le K} 
P_z\big( T(B(r-H_0))<T(\cV^c)\big)\\
\le & \sup_{\cV \subset\cS_0:\ |\cV|\le K}
C_d\exp\pare{-\kappa_d \pare{\frac{H_0^d}{|\cV|}}^{1/(d-1)}}\\
\le& C_d \exp\big(-\kappa_d \gamma^{1/(d-1)}\big).
\end{split}
\ee
Now, we define $\gamma$ large enough so that the right hand
side of \reff{explorer-1} is less than $1/e$.
Thus, we have an estimate valid for any explorer:
\[
P\big(\text{ an explorer crosses }\cS_0 \big) \le \frac{1}{e}.
\]
Also, each explorer having crossed $\cS_0$ is stopped upon
entering $B(2r-H_0)$. We call $\eta_1$ their configuration.
The key observation is that $|\eta_1|$ is bounded
by a Poisson variable $\cN_1$ with parameter $\lambda_0/e$.
Now, we define the second shell so that its width is
\[
H_1^d=\gamma \cN_1,\quad\text{and}\quad
\cS_1:=B(2r-H_0)\bs B(2r-H_1),
\]
and so forth. Thus, after considering $i\ge 1$ crossings,
we have a Poisson variable $\cN_i$ with parameter
$e^{-i}\lambda_0$ (with $\lambda_0= \epsilon r^d$)
bounding the number of explorers stopped upon
entering of $B(2r-H_0-\dots-H_{i-1})$, 
and a width $H_i^d=\gamma \cN_i$. This ensures that
any explorer has probability less than $1/e$ to
cross the $i$-th shell.

Now, note that the event that one of
the explorer of $\eta_0$ crosses $B(2r)\bs B(r)$ is 
contained in the event that
$\{H_0+\dots+H_L>r\}$ where $L=\inf\{k: H_k=0\}$.

Define $i^*$ to be the smallest integer so that 
$e^{-i^*} \epsilon r^d$
is less than 1. Our starting point is
\be{step-1}
\{\sum_{i<L} H_i>r\}\subset \{\sum_{i<i^*} H_i>\frac{r}{2}\}\cup
\{\sum_{i<i^*} H_i\le \frac{r}{2},
\ \sum_{i^*\le i<L} H_i \ge \frac{r}{2}\}.
\ee
Now, $\cN_{i^*}$ is bounded by a Poisson variable of parameter 1,
and we further divide the second event of \reff{step-1}
according to dimension.
\subsection{Dimension Two}
We write, for some small $\delta$
\be{key-d2}
\begin{split}
\{\sum_{i<L} H_i>r\}\subset& \{\sum_{i<i^*} H_i>\frac{r}{2}\}\cup
\{\sum_{i<i^*} H_i\le \frac{r}{2},
\ \sum_{i^*\le i<L} H_i \ge \frac{r}{2}\}\\
\subset& \{\sum_{i<i^*} \cN_i^{1/2} > 4 \rho \}\cup \{\cN_{i^*}> 
\frac{\delta r^2}{\log^2 r}\}\quad\text{with}\quad 4\rho:=
\frac{r}{2 \gamma^{1/2}}\\
& \cup \{\text{less than $\frac{\delta r^2}{\log^2 r}$ explorers
cross a shell of width } r/2\}.
\end{split}
\ee
We now proceed in estimating the three events separately. Note
that the event that less 
than $\frac{\delta r^2}{\log^2 r}$ explorers
have to cross a shell of width $r/2$ is dealt 
with Lemma~\ref{lem-levine}.
Also, the fact that $\cN_{i^*}$ 
is bounded by $\mathcal P(1)$ and \reff{ld-poisson} yield
\[
\begin{split}
P\big(\cN_{i^*}\ge \delta\frac{r^2}{\log^2 r}\big)\le&
\exp\big(- \delta\frac{r^2}{\log^2 r} \log(\frac{r^2}{e\log^2 r})\big)\\
\le & \exp\big(-2  \delta\frac{r^2}{\log r}(1+o(1))\big).
\end{split}
\]
We now deal with the deviation 
$\{\sum_{i<i^*} \cN_i^{1/2}>4\rho\}$.
A union bound allows us to treat this term
after we distinguish three regimes: (i) when $i$ is {\it small},
the deviation asks $\cN_i$ to be larger than $e^{-i} r^2$, and small
means for $i\le j^*$ with $\exp(j^*)\cdot\log(r)=1$, (ii) 
when $j^*\le i\le 2j^*$, we ask $\cN_i$ to be
larger than $r^2/(i^2\cdot \log^2 2)$,
and finally (iii) 
when $i$ is large the deviation asks $\cN_i$ 
to be larger than $r^2/(i\cdot i^*)$, and
we shall see that this gives the correct bound for $i\ge 2j^*$.
The first regime will fix the value for $\epsilon$.

More precisely, our first sum runs up to $j^*=\log(\log r)$, 
such that $\exp(j^*)=\log r$. We now use that
\[
\{\sum_{i< j^*}  \cN_i^{1/2}\ge \rho\}\subset \bigcup_{i< j^*} 
\{ \cN_i> \frac{e^{-i}}{(1-e^{-1/2})^2} \rho^2\}.
\]
For any $i\le j^*$, we have for $\epsilon$ small enough, and some
constant $\kappa$
\[
\begin{split}
P\big( \cN_i> \frac{e^{-i}}{(1-e^{-1/2})^2} \rho^2\big)\le&
\exp\big(-\frac{\rho^2}{(1-e^{-1/2})^2e^i}
\log\big(\frac{\rho^2}{e\epsilon(1-e^{-1/2})^2r^2}\big)\\
\le& \exp\big(-\kappa e^{-j^*} r^2\big)
=\exp\big(-\kappa\frac{r^2}{\log r}\big).
\end{split}
\]
Now, we consider case (ii). Note that 
\[
\{\sum_{j^*\le i\le 2j^*}  \cN_i^{1/2}\ge \rho\}
\subset \bigcup_{j^*\le i\le 2j^*}
\{ \cN_i> \frac{\rho^2}{i^2\cdot \log^2(2)}\}.
\]
Since $j^* = \log\log r$,
$e^{-i} i^2$ is smaller than $4(\log\log r)^2 / \log r$
for $j^* \leq i < 2j^*$.
Hence, by (3.1) there is a positive constant $c_0$
such that, for any such $i$,
$$
P\left(
{\cal N}_i > {\rho^2 \over i^2 \log^2 2}
\right) \leq \exp\left\{
- {\rho^2 \over i^2 \log^2 2}
\log\left( {c_0 \log r \over (\log\log r)^2} \right) \right\}
$$
and there is positive constant $\kappa$ such that
$$
P\left(
{\cal N}_i > {\rho^2 \over i^2 \log^2 2}
\right) \leq \exp\left\{
- \kappa {r^2 \over (\log\log r)^2}
\right\}.
$$
Finally, consider $2j^*\le i\le i^*$. First, note that
$$
\left\{\sum_{2j^* < i < i^*} {\cal N}_i^{1 / 2} > 2 \rho
\right\} \subset \bigcup_{2j^* < i < i^*} \left\{
{\cal N}_i^{1 / 2} > {\rho \over \sqrt{i^* \cdot i}}
\right\}.
$$
Secondly, there is a constant $C>1$, such that for
any $\epsilon$, and $r$ large enough
\[
i\cdot i^*\cdot e^{-i}\le \big(i\cdot e^{-i/4}\big)\times
\big(i^*\cdot e^{-j^*}\big)\times e^{-i/4}\le 
\frac{C\cdot e^{-i/4}}{e}.
\]
Indeed, we have used that $i^*$ is of order $2\log(r)$ whereas
$\exp(j^*)=\log(r)$.
Thus, for $2j^*\le i\le i^*$, there is a constant $\kappa>0$
such that for $r$ large enough
\[
\begin{split}
P\big( \cN_i>  \frac{\rho^2}{i^*\cdot i} \big)&\le
\exp\big(-\frac{\rho^2}{i\cdot i^*}\log(\frac{\rho^2}
{e\cdot i\cdot i^*\cdot e^{-i}r^2})\big)\\
&\le \exp\big(-\frac{\rho^2}{i^*}\times \frac{i/4+\log(C)}{i}\big)
\le \exp\big(-\kappa\frac{r^2}{\log(r)}\big).
\end{split}
\]
This concludes the proof of the Corollary in the case $d = 2$.

\subsection{Dimensions $d\ge 3$ }
We recall \reff{step-1}, and we write for some small $\delta$
with $2\rho:=r/(2\gamma^{1/d})$,
\be{key-d3}
\begin{split}
\{H_0+\dots+H_L>r\}\subset& \acc{\sum_{i<i^*} \cN_i^{1/d} >
2\rho}\cup \acc{\cN_{i^*}> \delta\frac{r^2}{
\log r}}\\
& \cup \acc{\text{less than $\frac{\delta r^2}{\log r}$ explorers
cross a shell of width } r/2}.
\end{split}
\ee
We now proceed in estimating the three events separately. 

Note that $\cN_{i^*}$ is bounded by $\mathcal P(1)$ so
\[
\begin{split}
P\big(\cN_{i^*}\ge \frac{\delta r^2}{\log r}\big)\le&
\exp\big(- \delta\frac{r^2}{\log r} \log(\frac{\delta r^2}{e\log r})\big)\\
\le & \exp\big(-2\delta r^2(1+o(1))\big).
\end{split}
\]
To deal with the deviation $\{\sum_{i<i^*} \cN_i^{1/d}>2\rho\}$, 
we use again a union bound, but here we only need 
to distinguish two regimes: (i) when $i$ is small,
the deviation asks $\cN_i$ to be larger than $c_0 e^{-i} r^d$, for some
fixed constant $c_0$ and thus $i$ small means that for some $\kappa>0$,
\[
P\big(\cN_i\ge c_0 e^{-i} r^d\big)\le \exp\big(-
c_0 e^{-i} r^d \log(c_0(e\epsilon)^{-1})\big)\ll \exp(-\kappa r^2).
\]
We define $j^*$ to be the largest integer $i$ such that
\[
c_0 r^{d-2} e^{-i} \log(c_0(e\epsilon)^{-1})\ge 1.
\]
Note that $j^*$ is of order $(d-2)\log r$, whereas
$i^*$ is of order $d\log r$. 
Now, for $i>j^*$, we use that for some constant $c_1>0$ such that
\[
\sum_{i=j^*}^{i^*}\frac{c_1}{i}\le 1.
\]
We use the estimate for $i>j^*$ and $i<i^*$, 
that for some constant $\kappa>0$, and $r$ large enough
\[
P\big(\cN_i^{1/d}\ge \frac{c_1}{i} \rho\big)\le \exp\big(-
c_1^d \frac{\rho^d}{(i^*)^d} \log(\frac{c_1^d \rho^{d}}
{e\epsilon(i^*)^d e^{-j^*}r^2)}\big)\le \exp(-\kappa r^2).
\]
\section{Proof of Lemma~\ref{lem-levine}}\label{sec-old}
For $C_d > 1$ and $\kappa_d> 1$ appearing in Lemma~\ref{lem-1}
we set
\be{nicolas}
    \gamma = \max\left(
        1, \left(
            2 \log C_d \over \kappa_d
        \right)^{d - 1}
    \right).
\ee
We divide $B(0, r)$ into shells
${\cal S}_0$, ${\cal S}_1$, \dots\ of widths
$H_0$, $H_1$, \dots\ 
We set $H_0 = h_0 = r / 4$,
and for $k \geq 1$ the width $H_k$ is random and depends
on the number $N_{k - 1}$ of explorers
settling in the previous shell.
$$
    H_k^d = \gamma N_{k - 1}.
$$
We denote by $L$
the first $k$ for which $\sum_{i < k} H_i > 3r / 4$
or $H_k < 1$, in which case $N_{k - 1} = 0$.
The shells are as follows. For $k <L$,
$$
    {\cal S}_{k}
    = B\Biggl(
        0, r - \sum_{i < k} H_i
    \Biggr) \setminus B\Biggl(
        0, r - \sum_{i \leq k} H_i
    \Biggr)
$$
For $r$ large enough, we have $\gamma |\eta| \leq (r / 4)^d = h_0^d$,
so that for all $k \geq 0$,
$$
 H_k^d \leq \gamma N_{k - 1} \leq \gamma |\eta|
    \leq \left(
        r \over 4
    \right)^d
    = h_0^d.
$$
Now, if $0 \in A(\eta)$,
then $\sum_{k < L} H_k$
has to be larger than $3r / 4$,
which implies, since $H_1 \leq H_0 = r / 4$,
that $\sum_{k = 2}^L H_k > r / 4$.
Also, for each $k < L$, $N_k$ explorers have to cross
the shells ${\cal S}_0$, ${\cal S}_1$, \dots,
${\cal S}_{k - 1}$.
Since $L < r$ we get by Lemma~1.1
that, writing $(n_k, k < l)$
for a generic family of $l$ positive integers
and writing $h_i$ for $(\gamma n_{i - 1})^{1 / d}$
if $i > 0$,
\begin{align}
    P\bigl( 0 \in A(\eta) \bigr) 
    &\leq \sum_{l < r} \sum_{(n_k, k < l)}
    {|\eta| \choose n_k, k <l}
    C_d^{\sum_{k < l} k n_k}
    \exp\left\{
        - \sum_{k < l} n_k \sum_{i < k} \kappa_d
        \left(
            h_i^d \over n_i
        \right)^{1 \over d - 1}
    \right\}  \nonumber\\
    &\qquad\qquad\qquad
    \times\1\left\{
        \sum_{k = 2}^l h_k > {r \over 4}
        \hbox{ and $h_k \leq h_0$ for all $k < l$}
    \right\} \label{julie}
\end{align}
Then, by the arithmetic-geometric inequality
and using~\eqref{nicolas},
it holds, for each $k > 0$ with $h_k \leq h_0$,
\begin{align*}
    {1 \over k} \sum_{i < k} \left(
        h_i^d \over n_i
    \right)^{1 \over d - 1}
 &\geq \left(
        \prod_{i < k} {h_i^d \over n_i}
    \right)^{1 \over k(d - 1)}
    = \left(
        {h_0^d \over n_{k - 1}}
        \gamma^{k - 1}
    \right)^{1 \over k(d - 1)}
    = \left(
        {h_0^d \over h_k^d}
        \gamma^k
    \right)^{1 \over k(d - 1)} \\
    &\geq \gamma^{1 \over d - 1}
    \geq {2 \log C_d \over \kappa_d}
\end{align*}
and as soon as $h_k \leq h_0$ for all $k < l$,
$$
    C_d^{\sum_{k < l} k n_k}
    \exp\left\{
        - \sum_{k < l} n_k \sum_{i < k} \kappa_d
        \left(
            h_i^d \over n_i
        \right)^{1 \over d - 1}
    \right\}
    \leq C_d^{-\sum_{k < l} k n_k}
    = C_d^{-{1 \over \gamma} \sum_{k < l} kh_{k + 1}^d}.
$$
By Hölder's inequality
$$
    \sum_{k = 2}^l h_k
    = \sum_{k = 1}^{l - 1} h_{k + 1}
    \leq \left(
        \sum_{k = 1}^{l - 1} k^{-{1 \over d - 1}}
    \right)^{d - 1 \over d}
    \left(
        \sum_{k = 1}^{l - 1} k h_{k + 1}^d
    \right)^{1 \over d}.
$$
Now, there is a positive constant $c_d$
such that in $d = 2$ 
$$
    \sum_{k < l} k h_{k + 1}^d
    \geq {
        \left(
            \sum_{k = 2}^l h_k
   \right)^2
    \over
        c_2 \log l
    },
$$
whereas if $d \geq 3$.
$$
    \sum_{k < l} k h_{k + 1}^d
    \geq {
        \left(
            \sum_{k = 2}^l h_k
        \right)^d
    \over
         c_d l^{d - 2}
    }.
$$
Hence, \eqref{julie} yields in $d=2$
$$
    P\bigl(
        0 \in A(\eta)
    \bigr) \leq r |\eta|^r r^{|\eta|} C_2^{
        - {r^2 \over 16 \gamma c_2 \log r}
    },
$$
and yields in $d\ge 3$
$$
    P\bigl(
        0 \in A(\eta)
    \bigr) \leq r |\eta|^r r^{|\eta|} C_d^{
        - {r^2 \over 4^d \gamma c_d}
    }.
$$
These bounds establish the required asymptotics
provided that $|\eta| \leq a r^2 / \log^2 r$
if $d = 2$, or $|\eta| \leq a r^2 / \log r$
if $d \geq 3$, for a small enough $a > 0$.

\section*{Acknowledgements}
We warmly thank Lionel Levine and Yuval Peres.

\end{document}